\documentclass{amsart}
\usepackage[dvips]{graphicx}
\usepackage{amscd}
\usepackage{amsmath}
\usepackage{amsxtra}
\usepackage{amsfonts}
\usepackage{amssymb}
\newtheorem{theorem}{Theorem}[section]
\newtheorem{corollary}[theorem]{Corollary}
\newtheorem{lemma}[theorem]{Lemma}
\newtheorem{proposition}[theorem]{Proposition}
\theoremstyle{definition}
\newtheorem{definition}[theorem]{Definition}
\newtheorem{remark}[theorem]{Remark}

\theoremstyle{remark}

\renewcommand{\theclaim}{\textup{\theclaim}}

\numberwithin{equation}{section}

\def\openone

{\mathchoice

{\hbox{\upshape \small1\kern-3.3pt\normalsize1}}

{\hbox{\upshape \small1\kern-3.3pt\normalsize1}}

{\hbox{\upshape \tiny1\kern-2.3pt\SMALL1}}

{\hbox{\upshape \Tiny1\kern-2pt\tiny1}}}

\makeatletter

\newbox\ipbox

\newcommand{\ip}[2]{\left\langle #1\,|\,#2\right\rangle}

\newcommand{\diracb}[1]{\left\langle #1\mathrel{\mathchoice

{\setbox\ipbox=\hbox{$\displaystyle \left\langle\mathstrut
#1\right.$}

\vrule height\ht\ipbox width0.25pt depth\dp\ipbox}

{\setbox\ipbox=\hbox{$\textstyle \left\langle\mathstrut
#1\right.$}

\vrule height\ht\ipbox width0.25pt depth\dp\ipbox}

{\setbox\ipbox=\hbox{$\scriptstyle \left\langle\mathstrut
#1\right.$}

\vrule height\ht\ipbox width0.25pt depth\dp\ipbox}

{\setbox\ipbox=\hbox{$\scriptscriptstyle \left\langle\mathstrut
#1\right.$}

\vrule height\ht\ipbox width0.25pt depth\dp\ipbox}

}\right. }

\newcommand{\dirack}[1]{\left. \mathrel{\mathchoice

{\setbox\ipbox=\hbox{$\displaystyle \left.\mathstrut
#1\right\rangle$}

\vrule height\ht\ipbox width0.25pt depth\dp\ipbox}

{\setbox\ipbox=\hbox{$\textstyle \left.\mathstrut
#1\right\rangle$}

\vrule height\ht\ipbox width0.25pt depth\dp\ipbox}

{\setbox\ipbox=\hbox{$\scriptstyle \left.\mathstrut
#1\right\rangle$}

\vrule height\ht\ipbox width0.25pt depth\dp\ipbox}

{\setbox\ipbox=\hbox{$\scriptscriptstyle \left.\mathstrut
#1\right\rangle$}

\vrule height\ht\ipbox width0.25pt depth\dp\ipbox}

} #1\right\rangle}

\newcommand{\ltwor}{L^{2}\left(\mathbb{R}\right)}

\newcommand{\ltworn}{L^{2}\left(\mathbb{R}^n\right)}

\newcommand{\ltwozn}{l^{2}\left(\mathbb{Z}^n\right)}

\newcommand{\Trace}{\operatorname*{Trace}}

\newcommand{\Tper}{\mathcal{T}_{per}}

\newcommand{\Rn}{\mathbb{R}^n}

\newcommand{\Zn}{\mathbb{Z}^n}

\input cyracc.def

\begin{document}
\title[Local trace function]{The Local Trace Function of Shift Invariant
Subspaces}
\author{Dorin Ervin Dutkay}
\address{Department of Mathematics\\
The University of Iowa\\
14 MacLean Hall\\
Iowa City, IA 52242-1419\\
U.S.A.}
\email{ddutkay@math.uiowa.edu}
\thanks{}

\begin{abstract}
We define the local trace function for subspaces of $\ltworn$ which are invariant under integer translation. Our trace
function contains the dimension function and the spectral function defined in \cite{BoRz} and completely
characterizes the given translation invariant subspace. It has properties such as positivity, additivity, monotony
and some form of
continuity. It behaves nicely under dilations and modulations. We use the local trace function to deduce, using
short and simple arguments, some
fundamental facts about wavelets such as the characterizing equations, the equality between the dimension function
and the multiplicity function and some new relations between scaling functions and wavelets.
\end{abstract}
\maketitle
\tableofcontents
\section{\label{Intro}Introduction}
Shift invariant spaces are closed subspaces of $\ltworn$ that are invariant under all integer translations
(also called shifts). They lie at the very heart of several areas such as the theory of wavelets, spline systems,
Gabor systems or approximation theory. Thus, a good understanding of the shift invariant spaces can prove itself
fruitful and give results in each of these areas. The local trace function is a new instrument for the
analysis of shift invariant spaces. It is associated to a shift invariant subspace of $\ltworn$ and a positive
operator on $\ltwozn$. And it confirms our expectations: when applied to a specific case it
gives results. Several fundamental facts about wavelets can be deduced quite easily with the aid of the
local trace
function ( we mention some of these facts: the equations that characterize
wavelets (remark \ref{rem5_2}), the equality between the dimension function and the multiplicity function
(remark \ref{rem5_3}) and some new equations that
relate multiscaling functions to multiwavelets (theorem \ref{th5_4} and corrolary \ref{cor5_5})). Another nice thing
about the local trace function is
that it includes, as special cases, the dimension function and the spectral function introduced by M.Bownik and Z.
Rzeszotnik in
\cite{BoRz} ( their paper was the main source of inspiration for us, several of the theorems and proofs
presented here are just extensions of the theorems and proofs from \cite{BoRz}). To get the dimension function,
just compute the local trace function using the operator $I$-the identity on $\ltwozn$. The spectral function
is just the local trace function associated to the projection $P_{\delta_0}$ onto the $0$-th component (see
proposition \ref{prop3_10}).
 \par
The local trace function is based on two main concepts: the range function introduced by Helson (\cite{H}) for
shift invariant subspaces of $\ltworn$ and the trace function for positive operators on $\ltwozn$. The definition
of the
local trace functions combines these concepts (see definition \ref{def3_1}). Another key fact is a simple
observation regarding the link between the trace function and normalized tight frames, namely that the trace
function can be computed as in its definition but using a normalized tight frame instead of an orthonormal basis
(see proposition \ref{prop21_2}). Moreover, the normalized tight frames are the only families of vectors that can be used to compute the trace
(see proposition \ref{prop21_4_1} and theorem \ref{th21_4_2}).
\par
We recall that a subset $\{e_i\,|\,i\in I\}$ of a Hilbert space $H$ is called
a frame with positive constants $A$ and $B$ if
$$A\|f\|^2\leq\sum_{i\in I}|\ip{f}{e_i}|^2\leq B\|f\|^2,\quad(f\in H).$$
If $A=B=1$ then it is called a normalized tight frame (or shortly NTF).
\par
In this paper we define the local trace function, investigate some
of its properties and apply it to wavelets to obtain, using only
short arguments, some known and some new results. In section
\ref{Therange} and section \ref{Thetrace} we recall some facts
about the range function and the trace function, respectively. We
establish the main properties of these functions, properties that
we need for section \ref{Thelocal} where we define and study the
main character of our paper: the local trace function. Many
properties shared by the range function and the trace function are
inherited by the local trace function: additivity (proposition
\ref{prop3_5} and \ref{prop3_6}), monotony (proposition
\ref{prop3_8}), nice behavior with respect to modulation and
dilation (proposition \ref{prop3_8} and \ref{prop3_9}).
\par
Also, there is a strong connection between our trace function and
the Gramian introduced and effectively used by A.Ron and Z.Shen in
\cite{RS1},\cite{RS2}, \cite{RS3} (see the remarks \ref{rem21_4_4}
and \ref{rem3_3_3}).
\par
Another fact that should be noticed is that the local trace function completely determines the shift invariant
subspace: two shift invariant subspaces are equal if and only if their local trace functions are equal
(proposition \ref{prop3_7}).
\par
Section \ref{convergence} contains some results about the behavior of the local trace function with respect to
limits. That is, when a sequence of shift invariant spaces has a limit (which is also shift invariant), in some
cases the local trace function of the limit space is the limit of the sequence of local trace functions. One of the
conditions is that the convergence is in the strong operator topology, but we put restrictions on the operator
(theorem \ref{th4_2}); the other theorem is a monotone convergence theorem (theorem \ref{th4_3}) and it works for
any positive operator on $\ltwozn$.
\par
In section \ref{wavelets}, we apply the local trace function to wavelets. Some very simple arguments lead to
important results: the equations that characterize wavelets (\ref{eq5_2_1}), (\ref{eq5_2_2}), the equality between
the multiplicity function and the dimension function (remark \ref{rem5_3}) and the relation between scaling
functions and
wavelets (relation (\ref{eq5_4_2})-(\ref{eq5_5_2})) are obtained just by writing some local trace function in two
ways.
\par
Before we engage in the analysis of the local trace function, we have to recall some definitions and theorems.
\par
The Fourier transform is given by
$$\widehat{f}(\xi)=\int_{\Rn}f(x)e^{-i\ip{x}{\xi}}\,dx,\quad(\xi\in\Rn).$$
If $V$ is closed subspace of a Hilbert space $H$ and $f\in H$ we denote by $P_V$ the projection onto $V$ and
by $P_f$ the operator defined by:
$$P_f(v)=\ip{v}{f}f,\quad(v\in H).$$
\begin{definition}\label{def1_1}
A closed subspace $V$ of $\ltwor$ is called shift invariant (or shortly SI) if
$$T_kV=V,\quad(k\in\Zn),$$
where $T_k$ is the translation by $k$ on $\ltworn$:
$$(T_kf)(\xi)=f(\xi-k),\quad(\xi\in\Rn,f\in\ltworn).$$
If $\mathcal{A}$ is a subset of $\ltworn$ then we denote by $S(\mathcal{A})$ the shift invariant space generated
by $\mathcal{A}$,
$$S(\mathcal{A})=\overline{\operatorname*{span}}\{T_k\varphi\,|\,k\in\Zn,\varphi\in\mathcal{A}\}.$$
\end{definition}
\begin{definition}\label{def1_2}
Let $V$ be a shift invariant subspace of $\ltworn$. A subset $\phi$ of $V$ is called a normalized tight frame
generator (or NTF generator) for $V$ if
$$\{T_k\varphi\,|\,k\in\Zn,\varphi\in\phi\}$$
is a NTF for $V$.
\par
We use also the notation $S(\varphi):=S(\{\varphi\})$. $\varphi$ is called a quasi-orthogonal generator for
$S(\varphi)$ if
$$\{T_k\varphi\,|\,k\in\Zn\}$$
is a NTF for $S(\varphi)$ and for all $\xi\in\Rn$,
$$\operatorname*{Per}|\widehat{\varphi}|^2(\xi):=\sum_{k\in\Zn}|\widehat{\varphi}|^2(\xi+2k\pi)\in\{0,1\}.$$
(actually, the second condition is a consequence of the first but we include it anyway).
\end{definition}
 Several proofs are available in the literature for the next theorem which guarantees the existence of NTF generators (see
 \cite{Bo1}, \cite{B}, \cite{BMM}).
\begin{theorem}\label{th1_3}
Suppose $V$ is a SI subspace of $\ltworn$. Then $V$ can be decomposed as an orthogonal sum
$$V=\bigoplus_{i\in\mathbb{N}}S(\varphi_i),$$
where $\varphi_i$ is a quasi-orthogonal generator of $S(\varphi_i)$. Moreover $\{\varphi_i\,|\, i\in\mathbb{N}\}$
is a NTF generator for $V$.
\end{theorem}
In fact, it is possible to derive a choice for the generators $\varphi_i$ in theorem \ref{th1_3} directly from the
Stone-spectral multiplicity theorem, as explained in \cite{BMM}.

 \section{\label{Therange}The range function}
\par
In this section we define the range function and state some of its properties. Some
preliminary notations are needed.
\par
The Hilbert space of square integrable vector functions $L^2\left([-\pi,\pi]^n,\ltwozn\right)$ consists of
all vector valued measurable functions $\phi: [-\pi,\pi]^n\rightarrow\ltwozn$ with the norm
$$\|\phi\|=\left(\int_{[-\pi,\pi]^n}\|\phi(\xi)\|_{\ltwozn}^2\,d\xi\right)^{1/2}<\infty.$$
The scalar product is given by
$$\ip{\phi}{\psi}:=\int_{[-\pi,\pi]^n}\ip{\phi(\xi)}{\psi(\xi)}_{\ltwozn}\,d\xi.$$
\par
The map $\mathcal{T}:\ltworn\rightarrow L^2\left([-\pi,\pi]^n,\ltwozn\right)$ defined for $f\in\ltworn$ by
$$\mathcal{T}f:[-\pi,\pi]^n\rightarrow\ltwozn,\quad\mathcal{T}f(\xi)=(\widehat{f}(\xi+2k\pi))_{k\in\Zn},\quad(\xi\in
[-\pi,\pi]^n)$$
is an isometric isomorphism (up to multiplication by $1/(2\pi)^{n/2}$) between $\ltworn$ and $L^2\left([-\pi,\pi]^n,\ltwozn\right)$.
\par
We will need some variations of these maps because they will make
the formulas nicer. We define $L_{per}^2\left(\Rn,\ltwozn\right)$
to be the space of measurable vector valued functions
$\phi:\Rn\rightarrow\ltwozn$ with the property that
$\phi|_{[-\pi,\pi]^n}$ belongs to
$L^2\left([-\pi,\pi]^n,\ltwozn\right)$ and they are periodic in
the following sense:
$$\phi(\xi+2k\pi)=\lambda(k)^*(\phi(\xi)),\quad(\xi\in\Rn),$$
where, for every $k\in\Zn$, $\lambda(k)$ is the shift by $k$ on $\ltwozn$ that is
$$(\lambda(k)\alpha)(l)=\alpha(l-k),\quad(l\in\Zn,\alpha\in\ltwozn).$$
The scalar product is defined by the same formula, the integral being taken over
$[-\pi,\pi]^n$.
\par
The map $\Tper:\ltworn\rightarrow L_{per}^2\left(\Rn,\ltwozn\right)$ is
defined by the same formula as $\mathcal{T}$,
$$\Tper f:\Rn\rightarrow\ltwozn,\,\Tper f(\xi)=(\widehat{f}(\xi+2k\pi))_{k\in\Zn},\quad (\xi\in\Rn)$$
but notice that $\xi$ is now in $\Rn$. Also, observe the periodicity property of $\Tper$:
$$\Tper f(\xi+2s\pi)=\lambda(s)^*(\Tper f(\xi)),\quad(\xi\in\Rn,s\in\Zn).$$

\begin{definition}\label{def2_0}
A range function is a measurable mapping
$$J:[-\pi,\pi]^n\rightarrow\{\mbox{ closed subspaces of }\ltwozn\}.$$
Measurable means weakly operator measurable, i.e., $\xi\mapsto\ip{P_{J(\xi)}a}{b}$ is measurable for any
choice of vectors $a,b\in\ltwozn$.
\par
A periodic range function is a measurable function
$$J_{per}:\Rn\rightarrow\{\mbox{ closed subspaces of }\ltwozn\},$$
with the periodicity property:
$$J_{per}(\xi+2k\pi)=\lambda(k)^*\left(J_{per}(\xi)\right),\quad(k\in\Zn,\xi\in\Rn).$$
\end{definition}
 Sometimes we will use the same letter to denote the subspace $J_{per}(\xi)$ and the projection onto
$J_{per}(\xi)$.
In terms of projections, the periodicity can be written as:
$$J_{per}(\xi+2k\pi)=\lambda(k)^*J_{per}(\xi)\lambda(k),\quad(k\in\Zn,\xi\in\Rn).$$
\par
 The next theorem due to Helson \cite{H} establishes the fundamental connection between shift
invariant spaces and the range function. The theorem appears in this form in \cite{Bo1}, proposition
1.5, the only modification needed is to work with the whole $\Rn$ instead of $[-\pi,\pi]^n$ and use
the periodicity.

\begin{theorem}\label{th2_1}
A closed subspace $V$ of $\ltworn$ is shift invariant if and only if
$$V=\{f\in\ltworn\,|\, \Tper f(\xi)\in J_{per}(\xi)\mbox{ for a.e. }\xi\in\mathbb{R}^n\},$$
for some measurable periodic range function $J_{per}$. The correspondence between $V$ and $J_{per}$ is
bijective under the convention that range functions are identified if they are equal a.e. Furthermore, if
$V=S(\mathcal{A})$ for some countable $\mathcal{A}\subset\ltworn$, then
$$J_{per}(\xi)=\overline{\operatorname*{span}}\{\Tper\varphi(\xi)\,|\,\varphi\in\mathcal{A}\},\quad \mbox{for
a.e. }\xi\in\Rn.$$
\end{theorem}

\begin{definition}\label{def2_1_1}
If $V$ is a SI subspace of $\ltworn$ then $J_{per}$ associated to $V$ as in theorem
$\ref{th2_1}$ is called the periodic range function of $V$.
\end{definition}

\par
The first elementary property of the range function that we will need is additivity.
The range function is also unitary in the sense that it preserves orthogonality of
subspaces. The precise formulation of these properties in given below.
\begin{proposition}\label{prop2_2}
Let $(V_i)_{i\in I}$ be a countable family of mutually orthogonal SI subspaces of
$\ltworn$ and denote by $J_{V_i}$ the periodic range function of $V_i$. If
$$V=\bigoplus_{i\in I}V_i$$
and $J_V$ is its periodic range function then
$$J_V(\xi)=\bigoplus_{i\in I}J_{V_i}(\xi),\quad\mbox{for a.e. }\xi\in\Rn,$$
where the sum is an orthogonal one.
\end{proposition}

\begin{proof}
Pick some countable $\phi_i\subset V_i$ such that $\{T_k\varphi\,|\, k\in\Zn,\varphi\in\phi_i\}$ spans $V_i$.
Then it is clear that $\{T_k\varphi\,|\, k\in\Zn,\varphi\in\phi\}$ spans $V$,
with the notation $\phi=\cup_{i\in I}\phi_i$.
Using Helson's theorem \ref{th2_1} we can determine the periodic range functions:
$$J_{V_i}(\xi)=\overline{\operatorname*{span}}\{\Tper\varphi(\xi)\,|\,\varphi\in\phi_i\},\quad(i\in
I),$$
$$J_{V}(\xi)=\overline{\operatorname*{span}}\{\Tper\varphi(\xi)\,|\,\varphi\in\phi\},$$
 for almost every point $\xi\in\Rn$. This shows that, if we check the orthogonality of the subspaces
$J_{V_i}(\xi)$, then we are done.
\par
To check that these subspaces of $\ltwozn$ are mutually orthogonal, take two arbitrary $i\neq j\in
I$ and $\varphi_1\in\phi_i$,$\varphi_2\in\phi_j$. Then, since we are dealing with SI
spaces,
$T_k\varphi_1$ is perpendicular to $\varphi_2$ for any choice of $k\in\Zn$. Rewriting
this in terms
of the Fourier transform, we obtain
$$0=\int_{\Rn}e^{-i\ip{k}{\xi}}\widehat{\varphi}_1(\xi)\overline{\widehat{\varphi}_2}(\xi)\,d\xi=
\int_{[-\pi,\pi]^n}e^{-i\ip{k}{\xi}}\ip{\mathcal{T}\varphi_1(\xi)}{\mathcal{T}\varphi_2(\xi)}\,d\xi.$$
For the second equality we applied a periodization. But this shows that all Fourier coefficients of
the map $\xi\mapsto\ip{\mathcal{T}\varphi_1(\xi)}{\mathcal{T}\varphi_2(\xi)}$ are 0 so
the map itself is 0 which implies that $\Tper\varphi_1(\xi)$ is perpendicular to
$\Tper\varphi_2(\xi)$ for a.e. $\xi\in[-\pi,\pi]^n$ and therefore, because of the
periodicity, for a.e.
$\xi\in\Rn$. Consequently, the subspaces are mutually orthogonal and the proposition is proved.
\end{proof}
\par
An easy consequence of the additivity of the range function is the following monotony
property:
\begin{proposition}\label{prop2_3}
Let $(V_j)_{j\in\mathbb{N}}$ be an increasing sequence of SI subspaces of $\ltworn$.
Denote by
$$V:=\overline{\bigcup_{j\in\mathbb{N}}V_j}.$$
Then for a.e. $\xi\in\Rn$, $(J_{V_j}(\xi))_{j\in\mathbb{N}}$ is increasing and
$$J_V(\xi)=\overline{\bigcup_{j\in\mathbb{N}}J_{V_j}(\xi)}.$$
\end{proposition}

\begin{proof}
Let $W_j$ be the orthogonal complement of $V_j$ in $V_{j+1}$, ($j\in\mathbb{N}$). Then $W_j$ are shift
invariant too, $V_{j+1}=V_{j}\oplus W_j$ for all $j\in\mathbb{N}$ and
$$V=V_0\oplus\bigoplus_{j\in\mathbb{N}}W_j,\quad
V_l=V_0\oplus\bigoplus_{j=0}^{l-1}W_j,\quad(l\in\mathbb{N}).$$
The proposition then follows from the additivity property stated in proposition \ref{prop2_2}.
\end{proof}

 The next proposition can also be found in \cite{H} and \cite{Bo1}. Again the only modification is that we
use the periodic extension of the range function.

\begin{proposition}\label{prop2_4}
Let $V$ be a SI subspace of $\ltworn$ and $J_{per}$ its periodic range function. Then
$$\Tper(P_Vf)(\xi)=J_{per}(\xi)(\Tper f(\xi)),\quad\mbox{ for a.e. }\xi\in\Rn.$$
\end{proposition}

\begin{definition}\label{def6_1}
A subset $\{e_i\,|\,i\in I\}$ of a Hilbert space $H$ is called a Bessel sequence with constant $B>0$ if
$$\sum_{i\in I}|\ip{e_i}{f}|^2\leq B\|f\|^2,\quad(f\in H).$$
\end{definition}

 The next theorem will be extensively used in this paper. Its proof can be found in
\cite{Bo1}.
\begin{theorem}\label{th2_5}
Let $V$ be a SI subspace of $\ltworn$, $J_{per}$ its periodic range function and $\phi$ a countable subset
of $V$.
$\{T_k\varphi\,|\,k\in\Zn,\varphi\in\phi\}$ is a frame with constants $A$ and $B$ for $V$ (Bessel family with constant
$B$) if and only if
$\{\Tper\varphi(\xi)\,|\,\varphi\in\phi\}$ is a frame with constants $A$ and $B$ for $J_{per}(\xi)$ (Bessel sequence
with constant $B$) for
almost every $\xi\in\Rn$.
\end{theorem}

\section{\label{Thetrace}The trace function}
In this section we recall the definition of the trace function and gather some of its properties, the main ones
and also some others that we will need in the sequel.
\begin{definition}\label{def21_1}
Let $H$ be a Hilbert space and $T$ a positive operator on $H$. Then the trace of the operator $T$ is the positive number (can
be also $\infty$) defined by
$$\Trace(T)=\sum_{i\in I}\ip{Te_i}{e_i},$$
where $\{e_i\,|\,i\in I\}$ is an orthonormal basis for $H$.
\end{definition}
\par
The next proposition shows that the trace is well defined, that is it doesn't depend on the choice of the orthonormal basis, and
moreover, it can be computed with the same formula using a normalized tight frame.
\begin{proposition}\label{prop21_2}
Let $H$ be a Hilbert space, $T$ a positive operator on $H$ and $\{f_j\,|\,j\in J\}$ a normalized tight frame for $H$. Then
$$\Trace(T)=\sum_{j\in J}\ip{Tf_j}{f_j}.$$
\end{proposition}

\begin{proof}
Let $\{e_i\,|\,i\in I\}$ be an orthonormal basis for $H$.
\begin{align*}
\Trace(T)&=\sum_{i\in I}\ip{Te_i}{e_i}=\sum_{i\in I}\ip{T^{1/2}e_i}{T^{1/2}e_i}=\sum_{i\in I}\|T^{1/2}e_i\|^2\\
&=\sum_{i\in I}\sum_{j\in J}|\ip{T^{1/2}e_i}{f_j}|^2=\sum_{j\in J}\sum_{i\in I}|\ip{T^{1/2}e_i}{f_j}|^2\\
&=\sum_{j\in J}\sum_{i\in I}|\ip{e_i}{T^{1/2}f_j}|^2=\sum_{j\in J}\|T^{1/2}f_j\|^2\\
&=\sum_{j\in J}\ip{Tf_j}{f_j}.
\end{align*}
\end{proof}
\par
The following proposition enumerates some of the elementary properties of the trace.
For a proof look in any basic book on operator theory (e.g. \cite{StZs}).
\begin{proposition}\label{prop21_3}
The trace has the following properties: for all $a,b$ positive operators and $\lambda\geq0$:
\begin{enumerate}
\item
$\Trace(a+b)=\Trace(a)+\Trace(b)$;
\item
$\Trace(\lambda a)=\lambda\Trace(a)$;
\item
If $a\leq b$ then $\Trace(a)\leq\Trace(b)$;
\item
$\Trace(v^*av)\leq\Trace(a)$ whenever $v$ is a partial isometry;
\item
$\Trace(u^*au)=\Trace(a)$ for all unitary $u$;
\item
$\Trace(x^*x)=\Trace(xx^*)$ for every operator $x$ on $H$;
\item
$\Trace(|x|)=\Trace(|x^*|)$ for every operator $x$ on $H$ ( $|x|:=(x^*x)^{1/2}$);
\item
$\Trace(P)=\operatorname*{dim}PH$ for any projection $P$ in $H$.
\end{enumerate}
\end{proposition}

\begin{proposition}\label{prop21_4}
Let $P$ be a projection and $\{e_i\,|\,i\in I\}$ a NTF for its range. Then, for any positive operator $T$, and any vector $f\in H$.
$$\Trace(TP)=\sum_{i\in I}\ip{Te_i}{e_i},$$
$$\Trace(P_fP)=\sum_{i\in I}|\ip{f}{e_i}|^2=\|Pf\|^2.$$
\end{proposition}
\begin{proof}
Let $\{f_j\,|\,j\in J\}$ be an orthonormal basis for the orthogonal complement of the range of $P$. Then
$\{e_i\,|\,i\in I\}\cup\{f_j\,|\,j\in J\}$ is a NTF for the entire Hilbert space. Using proposition
\ref{prop21_2}, the formula for $\Trace(TP)$ follows.
\par
Particularize the formula for $T=P_f$:
$$\Trace(P_fP)=\sum_{i\in I}\ip{\ip{e_i}{f}f}{e_i}=\sum_{iu\in I}|\ip{f}{e_i}|^2=\|Pf\|^2.$$
\end{proof}

In fact even more is true: the equations of proposition
\ref{prop21_4} characterize the normalized tight frames for the
range of $P$. This is described in the next statement.
\begin{proposition}\label{prop21_4_1}
Let $H$ be a Hilbert space, $H_0$ a closed subspace and $\{e_i\,|\,i\in I\}$ a family of vectors from $H$. The following affirmations are equivalent:
\begin{enumerate}
\item
$\{e_i\,|\,i\in I\}$ is a NTF for $H_0$;
\item
For every positive operator $T$ on $H$,
$$\Trace(TP_{H_0})=\sum_{i\in I}\ip{Te_i}{e_i};$$
\item
For every vector $v\in H$,
$$\Trace(P_vP_{H_0})=\sum_{i\in I}|\ip{v}{e_i}|^2.$$
\end{enumerate}
(Note that we do not require in (ii) and (iii) that the $e_i$'s be in $H_0$. This will follow from the formulas).
\end{proposition}

\begin{proof}
(i) implies (ii) and (ii) implies (iii) according to proposition \ref{prop21_4_1} and its proof. So we only have to worry about the implication from
(iii) to (i). The hypotheses imply that, if $v\perp H_0$, then $P_vP_{H_0}=0$ so $v$ is perpendicular to $e_i$ for every $i\in I$.
\par
Also, for $v\in H_0$,
$$\sum_{i\in I}|\ip{v}{e_i}|^2=\|v\|^2.$$
This shows that all $e_i$'s are in $H_0$ and they form a NTF for it.
\end{proof}

\par
We can weaken the condition (iii) in proposition \ref{prop21_4_1}. The family of vectors is a NTF for the subspace
if it satisfies the NTF condition just for some special vectors as shown below:
\begin{theorem}\label{th21_4_2}
Let $H$ be a Hilbert space and $\{\delta_k\,|\,k\in K\}$ a total subset of $H$ (i.e. its closed linear span is $H$). Let $H_0$ be a closed subspace for
$H$ and
$\{e_i\,|\,i\in I\}$ a family of vectors in $H$. The following affirmations are equivalent:
\begin{enumerate}
\item
$\{e_i\,|\,i\in I\}$ is a NTF for $H_0$;
\item
For every $r\neq s\in K$ and $\lambda\in\{0,1,i=\sqrt{-1}\}$
\begin{equation}\label{eq21_4_2_1}
\sum_{i\in I}|\ip{\delta_r+\lambda\delta_s}{e_i}|^2=\|P_{H_0}(\delta_r+\lambda\delta_s)\|^2.
\end{equation}
\end{enumerate}
\end{theorem}

\begin{proof}
(i) implies (ii) clearly.
\par
Assume (ii).
Take $\{f_j\,|\,j\in J\}$ a NTF for $H_0$. Then
\begin{equation}\label{eq21_4_2_2}
\|P_{H_0}(v)\|^2=\sum_{j\in J}|\ip{v}{f_j}|^2,\quad(v\in H).
\end{equation}
Take $\lambda=0$ and use (\ref{eq21_4_2_1}) and (\ref{eq21_4_2_2}):
\begin{equation}\label{eq21_4_2_3}
\sum_{i\in I}|\ip{\delta_r}{e_i}|^2=\sum_{j\in J}|\ip{\delta_r}{f_j}|^2,\quad(v\in H).
\end{equation}
Now take $r\neq s$ and $\lambda\in\{1,i\}$. From (\ref{eq21_4_2_1}) and (\ref{eq21_4_2_2}):
$$\sum_{i\in I}(|\ip{\delta_r}{e_i}|^2+\ip{\delta_r}{e_i}\overline{\lambda}\overline{\ip{\delta_s}{e_i}}+
\overline{\ip{\delta_r}{e_i}}\lambda\ip{\delta_s}{e_i}+|\lambda|^2|\ip{\delta_s}{e_i}|^2)=$$
$$\sum_{j\in J}(|\ip{\delta_r}{f_j}|^2+\ip{\delta_r}{f_j}\overline{\lambda}\overline{\ip{\delta_s}{f_j}}+
\overline{\ip{\delta_r}{f_j}}\lambda\ip{\delta_s}{f_j}+|\lambda|^2|\ip{\delta_s}{f_j}|^2).$$
With (\ref{eq21_4_2_3}) we can reduce this to
$$\overline{\lambda}\sum_{i\in I}\ip{\delta_r}{e_i}\overline{\ip{\delta_s}{e_i}}+\lambda\sum_{i\in
I}\overline{\ip{\delta_r}{e_i}}\ip{\delta_s}{e_i}=$$
$$=\overline{\lambda}\sum_{j\in J}\ip{\delta_r}{f_j}\overline{\ip{\delta_s}{f_j}}+\lambda\sum_{j\in
J}\overline{\ip{\delta_r}{f_j}}\ip{\delta_s}{f_j}$$
Now take $\lambda=1$ and $\lambda=i$ and the two resulting equation will yield:
\begin{equation}\label{eq21_4_2_4}
\sum_{i\in I}\ip{\delta_r}{e_i}\overline{\ip{\delta_s}{e_i}}=\sum_{j\in J}\ip{\delta_r}{f_j}\overline{\ip{\delta_s}{f_j}}
\end{equation}
Next we prove that
$$\sum_{i\in I}|\ip{v}{e_i}|^2=\sum_{j\in J}|\ip{v}{f_j}|^2$$
holds for all $v\in S$ where
$$S:=\{\sum_{k\in K_0}v_k\delta_k\,|\,K_0\subset K\mbox{ finite }\}.$$
For this take an arbitrary $v=\sum_{k\in K_0}v_k\delta_k$, with $K_0$ finite. Then
\begin{align*}
\sum_{i\in I}|\ip{v}{e_i}|^2&=\sum_{i\in I}|\sum_{k\in K_0}v_k\ip{\delta_k}{e_i}|^2\\
&=\sum_{i\in I}\sum_{k,k'\in K_0}v_k\overline{v}_{k'}\ip{\delta_k}{e_i}\overline{\ip{\delta_{k'}}{e_i}}\\
&=\sum_{k,k'\in K_0}v_k\overline{v}_{k'}\sum_{i\in I}\ip{\delta_k}{e_i}\overline{\ip{\delta_{k'}}{e_i}}\\
&=\sum_{k,k'\in K_0}v_k\overline{v}_{k'}\sum_{j\in J}\ip{\delta_k}{f_j}\overline{\ip{\delta_{k'}}{f_j}}\quad(\mbox{with (\ref{eq21_4_2_4})})\\
&=\sum_{j\in J}|\ip{v}{f_j}|^2.
\end{align*}
To prove the relation for arbitrary $v\in H$, define $\tilde{T}_1:S\rightarrow l^2(I)$, $T_2:H\rightarrow l^2(J)$ by
$$\tilde{T}_1v=(\ip{v}{e_i})_{i\in I},\quad(v\in S),\quad T_2v=(\ip{v}{f_j})_{j\in J}.$$
Then $\tilde{T}_1$ is a well defined linear operator, $\|T_2v\|^2=\|P_{H_0}v\|^2$ for $v\in H$ and
$\|\tilde{T}_1v\|=\|T_2v\|$ for $v\in S$. This shows that $\|\tilde{T}_1v\|^2\leq\|v\|^2$ for $v\in S$.  Hence, as
$S$ is a dense subspace of $H$, we can extend $\tilde{T}_1$ to a linear operator $T_1$ on $H$ with
$$\|T_1v\|^2\leq\|v\|^2,\quad\mbox{for all }v\in H.$$
We claim that $T_1v=(\ip{v}{e_i})_{i\in I}$ for any $v\in H$. Indeed, $v$ can be approximated by vectors $v_n$ in $S$. Then,
for each $i\in I$,
$$(T_1v)_i=\lim_{n\rightarrow\infty}(T_1v_n)_i=\lim_{n\rightarrow\infty}\ip{v_n}{e_i}=\ip{v}{e_i}.$$
Also, taking the limit,
$$\|T_1v\|^2=\|T_2v\|^2=\|P_{H_0}v\|^2,\quad(v\in H),$$
and this implies (i).
\end{proof}

\begin{theorem}\label{th21_4_3}
Let $H$ be a Hilbert space $H_0$ a closed subspace , $\{\delta_k\,|\,k\in K\}$ a total set for $H$, and $\{f_j\,|\,j\in J\}$ a NTF for $H_0$. Let
$\{e_i\,|\,i\in I\}$ be a family of vectors in
$H$. The following affirmations are equivalent:
\begin{enumerate}
\item
$\{e_i\,|\,i\in I\}$ is a NTF for $H_0$;
\item
For all $r,s\in H_0$,
$$\sum_{i\in I}\ip{\delta_r}{e_i}\overline{\ip{\delta_s}{e_i}}=\sum_{j\in J}\ip{\delta_r}{f_j}\overline{\ip{\delta_s}{f_j}}.$$
\end{enumerate}
\end{theorem}
\begin{proof}
Just examine the proof of theorem \ref{th21_4_2}.
\end{proof}

\begin{remark}\label{rem21_4_4}
Theorem \ref{th21_4_3} shows that $\{e_i\,|\, i\in I\}$ is a NTF
for $H_0$ if and only if it has the same Gramian as
$\{f_j\,|\,j\in J\}$. We recall briefly this notion. The Gramian
was introduced by A.Ron and Z.Shen in a series of papers and they
used it to analize the structure of shift invariant spaces
(\cite{RS1}, \cite{RS2}, \cite{RS3}).
\par
For a given countable family of vectors $\{e_i\,|\,i\in I\}$ define the operator $K:H\rightarrow l^2(I)$,
initially on sequences $f=(f_i)_{i\in I}$ with compact support, by
\begin{equation}\label{eq6_1}
K(f)=\sum_{i\in I}f_ie_i.
\end{equation}
If $K$ extends to a bounded operator (this is the case when $\{e_i\,|\,i\in I\}$ forms a Bessel family), then its
adjoint $K^*:l^2(I)\rightarrow H$ is given by
\begin{equation}\label{eq6_2}
K^*(v)=(\ip{v}{e_i})_{i\in I},\quad(v\in\ltwozn).
\end{equation}
\par
Suppose $\{e_i\,|\,i\in I\}$ is a subset of $H$ and $K$ is defined
as in (\ref{eq6_1}). The Gramian of the system $\{e_i\,|\,i\in
I\}$ is $G:l^2(I)\rightarrow l^2(I)$ defined by $G=K^*K$. The dual
Gramian of the system $\{e_i\,|\,i\in I\}$ is
$\tilde{G}:H\rightarrow H$ defined by $\tilde{G}=KK^*$.
\par
Note that
\begin{equation}\label{eq6_2_1}
\ip{\tilde{G}f}{g}=\ip{K^*f}{K^*g}=\sum_{i\in I}\ip{f}{e_i}\overline{\ip{g}{e_i}}=:B(f,g).
\end{equation}
The trace function involves expressions of the form
$$\sum_{i\in I}|\ip{f}{e_i}|^2=B(f,f)=\ip{\tilde{G}f}{f}.$$
We can recuperate the dual Gramian if these expressions are given,
because by polarization
$$B(f,g)=\frac{1}{4}\sum_{k=0}^3i^kB(f+i^kg,f+i^kg).$$
$B$, as a sesquilinear form gives rise to an operator and this
operator has to be the dual Gramian $\tilde{G}$. Another fact that
is worth noticing is that, if $\{e_i\,|\,i\in I\}$ is a NTF for
some subspace $H_0$ of $\ltwozn$, then the dual Gramian
$\tilde{G}$ is the projection onto $H_0$, $P_{H_0}$. This is
because $P_{H_0}$ verifies the equation (\ref{eq6_2_1}).
\end{remark}

\begin{proposition}\label{prop21_5}
Let $(P_j)_{j\in J}$ be a family of mutually orthogonal projections in a Hilbert space $H$. Then
$$\Trace(T\sum_{j\in J}P_j)=\sum_{j\in J}\Trace(TP_j),$$
for any positive operator $T$ on $H$.
\end{proposition}

\begin{proof}
Consider $\{e_i\,|\,i\in I_j\}$ an orthonormal basis for the range of $P_j$, ($j\in J$). Then, $\{e_i\,|\,i\in\cup
I_j\}$ is an orthonormal basis for the range of $\sum_{j\in J}P_j$. The additivity property is now clear, if we
use proposition \ref{prop21_4}.
\end{proof}
\par
The next proposition is an easy consequence of the additivity of the trace mentioned in proposition
\ref{prop21_5}.
\begin{proposition}\label{prop21_6}
Let $(P_j)_{j\in\mathbb{N}}$ be an increasing sequence of projections in some Hilbert space $H$,
$P=\sup_{j\in\mathbb{N}}P_j$ and $T$ a
positive operator on $H$. Then \\
$(\Trace(TP_j))_{j\in\mathbb{N}}$ increases to
$\Trace(TP)$.
\end{proposition}

\begin{proposition}\label{prop21_6_1}
Let $T$ be a positive operator on a Hilbert space $H$, $P$ a projection on $H$ and $U$
a unitary on $H$. Then
$$\Trace(UTPU^*)=\Trace(TP).$$
\end{proposition}
\begin{proof}
Let $\{e_i\,|\,i\in I\}$ be an orthonormal basis for $H$ such that $\{e_i\,|\, i\in I_0\}$
is an orthonormal basis for $UPH$, with $I_0\subset I$. Then $\{U^*e_i\,|\,i\in I_0\}$
is an orthonormal basis for $PH$ and $PU^*e_i=0$ for $i\in I\setminus I_0$; so
$$\Trace(UTPU^*)=\sum_{i\in I_0}\ip{UTPU^*e_i}{e_i}=\sum_{i\in
I_0}\ip{TPU^*e_i}{U^*e_i}=\Trace(TP).$$
\end{proof}
\begin{lemma}\label{lem21_8_1}
If $H$ is a Hilbert space, $H_0$ a closed subspace and $f\in H$, then
$\Trace(P_fP_{H_0})\leq\|f\|^2$ with equality if and only if $f\in H_0$.
\end{lemma}

\begin{proof}
Let $\{e_i\,|\, i\in I\}$ be an orthonormal basis for $H$ such that $\{e_i\,|\, i\in I_0\}$ is an orthonormal basis
for $H_0$ with $I_0\subset I$. Then
\begin{align*}
\|f\|^2&=\sum_{i\in I_0}|\ip{f}{e_i}|^2+\sum_{i\in I\setminus I_0}|\ip{f}{e_i}|^2\\
&=\sum_{i\in I_0}\ip{P_fe_i}{e_i}+\sum_{i\in I\setminus I_0}|\ip{f}{e_i}|^2\\
&=\Trace (P_fP_{H_0})+\sum_{i\in I\setminus I_0}|\ip{f}{e_i}|^2
\end{align*}
Thus the inequality holds and $\Trace(P_fP_{H_0})=\|f\|^2$ iff $\ip{f}{e_i}=0$ for all $i\in I\setminus I_0$
which is equivalent to $f\in H_0$.
\end{proof}

\begin{lemma}\label{lem21_9_1}
Let $H$ be a Hilbert space, $H_1$ and $H_2$ two closed subspaces. Then $H_1\subset H_2$ if and only if $\Trace(P_fP_{H_1})\leq
\Trace(P_fP_{H_2})$ for all $f\in H$.
\end{lemma}

\begin{proof}
If $H_1$ is contained in $H_2$ then, using proposition \ref{prop21_5}, the inequality
between the traces is immediate.
\par
For the converse, assume $f\in H_1$. Then, by lemma \ref{lem21_8_1}, $\Trace(P_fP_{H_1})=\|f\|^2$ so
$\|f\|^2\leq\Trace(P_fP_{H_2})$. But then, with the same lemma, $\Trace(P_fP_{H_2})=\|f\|^2$ which
implies that $f\in H_2$. As $f$ was arbitrary, the inclusion is proved.
\end{proof}

\section{\label{Thelocal}The local trace function}

\begin{definition}\label{def3_1}
Let $V$ be a SI subspace of $\ltworn$, $T$ a positive operator on $\ltwozn$ and let $J_{per}$ be the range function
associated to $V$. We define the local trace function associated to $V$ and $T$ as the map from $\mathbb{R}^n$ to
$[0,\infty]$ given by the formula
$$\tau_{V,T}(\xi)=\Trace\left(TJ_{per}(\xi)\right),\quad(\xi\in\mathbb{R}).$$
We define the restricted local trace function associated to $V$ and a vector $f$ in $\ltwozn$ by
$$\tau_{V,f}(\xi)=\Trace\left(P_fJ_{per}(\xi)\right)(=\tau_{V,P_f}(\xi)),\quad(\xi\in\mathbb{R}^n),$$
where $P_f$ is the operator on $\ltwozn$ defined by $P_f(v)=\ip{v}{f}f$.
\end{definition}

\begin{proposition}\label{prop3_1}
For all $f\in\ltwozn$,
$$\tau_{V,f}(\xi)=\|J_{per}(\xi)(f)\|^2,\quad(\mbox{ for a.e. }\xi\in\mathbb{R}^n).$$
\end{proposition}
\begin{proof}
Use the proposition \ref{prop21_4}.
\end{proof}

\begin{theorem}\label{th3_3}
Let $V$ be a SI subspace of $\ltworn$ and $\phi\subset V$ a NTF generator for $V$. Then for every positive operator $T$
on $\ltwozn$ and any $f\in\ltwozn$,
\begin{equation}\label{eq3_3_1}
\tau_{V,T}(\xi)=\sum_{\varphi\in\phi}\ip{T\Tper\varphi(\xi)}{\Tper\varphi(\xi)},\quad(\mbox{
for a.e. }\xi\in\mathbb{R}^n);
\end{equation}
\begin{equation}\label{eq3_3_2}
\tau_{V,f}(\xi)=\sum_{\varphi\in\phi}|\ip{f}{\Tper\varphi(\xi)}|^2,\quad \mbox{ for a.e. }\xi\in\mathbb{R}^n).
\end{equation}
\end{theorem}

\begin{proof}
According to theorem \ref{th2_5}, $\{\Tper\varphi(\xi)\,|\,\varphi\in\phi\}$ is a NTF for $J_{per}(\xi)$ for
a.e.
$\xi\in\mathbb{R}^n$.
With the proposition \ref{prop21_4} we obtain (\ref{eq3_3_1}). Take $T=P_f$ and
(\ref{eq3_3_1})
becomes (\ref{eq3_3_2}).
\end{proof}

\begin{theorem}\label{th3_3_1}
Let $V$ be a SI subspace of $\ltworn$, $J_{per}$ its periodic range function and $\phi$ a countable subset of $\ltworn$.
Then following affirmations are
equivalent:
\begin{enumerate}
\item
$\phi\subset V$ and $\phi$ is a NTF generator for $V$;
\item
For every $f\in\ltwozn$
\begin{equation}\label{eq3_3_1_1}
\sum_{\varphi\in\phi}|\ip{f}{\Tper\varphi(\xi)}|^2=\|J_{per}(\xi)(f)\|^2,\quad\mbox{for a.e. }\xi\in\Rn
\end{equation}
\item
For every $0\neq l\in\Zn$ and $\alpha\in\{0,1,i\}$,
\begin{equation}\label{eq3_3_1_2}
\sum_{\varphi\in\phi}|\widehat{\varphi}(\xi)+\overline{\alpha}\widehat{\varphi}(\xi+2l\pi)|^2=\|J_{per}(\xi)(\delta_0+\alpha\delta_l)\|^2,\quad\mbox{for
a.e. }\xi\in\Rn.
\end{equation}
\end{enumerate}
\end{theorem}

\begin{proof}
(i) implies (ii) just as an application of theorem \ref{th3_3} and proposition \ref{prop3_1}. (iii) is just a particular
case of (ii), namely $f=\delta_0+\alpha\delta_l$. So assume (iii) holds. Then (\ref{eq3_3_1_2}) can be rewritten as:
$$\sum_{\varphi\in\phi}{\ip{\delta_0+\alpha\delta_l}{\Tper(\xi)}|^2=\|J_{per}(\delta_0+\alpha\delta_l)}\|^2.$$
Now take $r\neq s\in\Zn$ and apply this equation to $l=s-r$, $\xi=\xi+2\pi r$. A short computation, that uses the
periodicity of $\Tper$ and $J_{per}$ and the fact that $\lambda(r)$ is unitary, will lead to
$$\sum_{\varphi\in\phi}|\ip{\delta_r+\alpha\delta_s}{\Tper(\xi)}|^2=\|J_{per}(\xi)(\delta_r+\alpha\delta_s)\|^2,\quad\mbox{
for a.e. }\xi\in\Rn.$$
Now we can use theorem \ref{th21_4_2} to conclude that, for a.e. $\xi\in\Rn$,
$$\{\Tper\varphi(\xi)\,|\,\varphi\in\phi\}$$
is a NTF for $J_{per}(\xi)$. And, with theorem \ref{th2_5}, (i) is obtained.
\end{proof}
Applying this theorem to $V=\ltworn$ we deduce the next corollary
which can be found also in \cite{RS1}:

\begin{corollary}\label{cor3_3_1}
A countable subset $\varphi$ of $\ltworn$ is a NTF generator for $\ltworn$ if and only if the following equations hold for
a.e. $\xi\in\Rn$:
$$\sum_{\varphi\in\phi}|\widehat{\varphi}(\xi)|^2=1,$$
$$\sum_{\varphi\in\phi}\widehat{\varphi}(\xi)\overline{\widehat{\varphi}}(\xi+2l\pi)=0,(l\in\Zn,l\neq 0).$$
\end{corollary}

\begin{theorem}\label{th3_3_2}
Let $V$ be a SI subspace of $\ltworn$, $\phi_1$ a NTF generator for $V$ and $\phi_2$ a countable family of vectors from
$\ltworn$. The following affirmations are equivalent:
\begin{enumerate}
\item
$\phi_2\subset V$ and $\phi_2$ is a NTF generator for $V$;
\item
For every $l\in\Zn$,
$$\sum_{\varphi\in\phi_2}\widehat{\varphi}(\xi)\overline{\widehat{\varphi}}(\xi+2l\pi)=
\sum_{\varphi\in\phi_1}\widehat{\varphi}(\xi)\overline{\widehat{\varphi}}(\xi+2l\pi),\quad\mbox{for a.e. }\xi\in\Rn.$$
\end{enumerate}
\end{theorem}

\begin{proof}
By theorem \ref{th2_5}, (i) is equivalent to
\begin{equation}\label{eq3_3_2_1}
\{\Tper\varphi(\xi)\,|\,\varphi\in\phi_2\}\mbox{ is a NTF for }J_{per}(\xi)\mbox{ for a.e. }\xi\in\Rn.
\end{equation}
But we know that $\{\Tper\varphi(\xi)\,|\,\varphi\in\phi_1\}$ is a NTF for $J_{per}(\xi)$ (again, by theorem
\ref{th2_5}). Therefore, using theorem \ref{th21_4_3}, (\ref{eq3_3_2_1}) is equivalent to:
\par
For every $r,s\in\Zn$ and a.e. $\xi\in\Rn$,
$$\sum_{\varphi\in\phi_2}\ip{\delta_r}{\Tper\varphi(\xi)}\overline{\ip{\delta_s}{\Tper\varphi(\xi)}}
=\sum_{\varphi\in\phi_1}\ip{\delta_r}{\Tper\varphi(\xi)}\overline{\ip{\delta_s}{\Tper\varphi(\xi)}},$$
which is exactly
$$\sum_{\varphi\in\phi_2}\widehat{\varphi}(\xi+2r\pi)\overline{\widehat{\varphi}}(\xi+2s\pi)=
\sum_{\varphi\in\phi_1}\widehat{\varphi}(\xi+2r\pi)\overline{\widehat{\varphi}}(\xi+2s\pi),\quad(r,s\in\Zn,\xi\in\Rn).$$
This implies (ii) and (ii) implies this, because one can take $l=s-r$, $\xi=\xi+2r\pi$.
\end{proof}

\begin{remark}\label{rem3_3_3}
Theorem \ref{th3_3_2} shows that $\phi_2$ is a NTF generator for
$V$ iff it has the same dual Gramian as $\phi_1$. Recall, that the
dual Gramian of a countable subset $\phi$ of $\ltworn$ is defined
as the function which assigns to each $\xi\in\Rn$ the dual Gramian
of the set $\{\Tper\varphi(\xi)\,|\,\varphi\in\phi\}$ (see remark
\ref{rem21_4_4} and \cite{RS1}, \cite{RS2}, \cite{RS3}, \cite{Bo1}
for details). The dual Gramian satisfies the equations:
$$\ip{\tilde{G}(\xi)\delta_r}{\delta_s}=\sum_{\varphi\in\phi}\widehat{\varphi}(\xi+2\pi r)\overline{\widehat{\varphi}}(\xi+2\pi s).$$
\end{remark}

\begin{proposition}\label{prop3_4}
{\bf [Periodicity]}
Let $V$ be a SI subspace, $T$ a positive operator on $\ltwozn$, $f\in\ltwozn$. Then, for $k\in\mathbb{Z}^n$
\begin{equation}\label{eq3_4_1}
\tau_{V,T}(\xi+2k\pi)=\tau_{V,\lambda(k)T\lambda(k)^*}(\xi),\quad
(\mbox{ for a.e. }\xi\in\mathbb{R}^n);
\end{equation}
\begin{equation}\label{eq3_4_2}
\tau_{V,f}(\xi+2k\pi)=\tau_{V,\lambda(k)f}(\xi),\quad (\mbox{ for a.e. }\in\mathbb{R}^n),
\end{equation}
where $\lambda(k)$ is the unitary operator on $\ltwozn$ defined by
$(\lambda(k)\alpha)(l)=\alpha(l-k)$ for all
$l\in\mathbb{Z}^n$, $\alpha\in\ltwozn$.
\end{proposition}

\begin{proof}
The periodicity of the local trace function is a consequence of the periodicity of the range function. Indeed,
by definition \ref{def2_0}, we know that for a.e. $\xi\in\mathbb{R}^n$ and every
$k\in\mathbb{Z}^n$,
$$J_{per}(\xi+2k\pi)=\lambda(k)^*J_{per}(\xi)\lambda(k).$$
Apply this in the definition of the local trace function and use
proposition \ref{prop21_6_1}:
$$\tau_{V,T}(\xi+2k\pi)=\Trace(TJ_{per}(\xi+2k\pi))=\Trace(T\lambda(k)^*J_{per}(\xi)\lambda(k))$$
$$=\Trace(\lambda(k)T\lambda(k)^*J_{per}(\xi))=\tau_{V,\lambda(k)T\lambda(k)^*}(\xi).$$
Since $\lambda(k)P_f\lambda(k)^*=P_{\lambda(k)f}$, (\ref{eq3_4_2}) follows from (\ref{eq3_4_1}).
\end{proof}

\begin{proposition}\label{prop3_5}
Let $V$ be a SI space, $T,S$ positive operators on $\ltwozn$, $f\in\ltwozn$
\begin{enumerate}
\item
$0\leq\tau_{V,T}(\xi)\leq\infty$, $0\leq\tau_{V,f}\leq\|f\|^2$;
\item
$\tau_{V,T+S}=\tau_{V,T}+\tau_{V,S}$;
\item
$\tau_{V,\lambda T}=\lambda\tau_{V,T},\quad (\lambda>0)$;
\item
$\tau_{V,\lambda f}=|\lambda|^2\tau_{V,f},\quad(\lambda\in\mathbb{C}).$
\end{enumerate}
\end{proposition}
\begin{proof}
Everything follows from section \ref{Thetrace}.
\end{proof}
\begin{proposition}\label{prop3_6}
{\bf [Additivity]}
Suppose $(V_i)_{i\in I}$ are mutually orthogonal SI subspaces ($I$ countable) and let $V=\oplus_{i\in I}V_i$. Then,
for every positive operator $T$ on $\ltwozn$ and every $f\in\ltwozn$:
\begin{equation}\label{eq3_6_1}
\tau_{V,T}=\sum_{i\in I}\tau_{V_i,T},\quad\mbox{a.e. on }\mathbb{R}^n;
\end{equation}
\begin{equation}\label{eq3_6_2}
\tau_{V,f}=\sum_{i\in I}\tau_{V_i,f},\quad\mbox{a.e. on }\mathbb{R}^n.
\end{equation}
\end{proposition}

\begin{proof}
Let $J_i$ be the periodic range function of $V_i$, ($i\in I$) and $J$ the periodic
range function of $V$. The range function
is
additive (proposition \ref{prop2_2}) so
$$J(\xi)=\sum_{i\in I}J_i(\xi),\quad\mbox{ a.e. on }\mathbb{R}^n.$$
Also the trace has additive properties (proposition \ref{prop21_5}) and these imply the
additivity of the local trace function
expressed in (\ref{eq3_6_1}). Again (\ref{eq3_6_2}) is just a particular case of (\ref{eq3_6_1}).
\end{proof}

\begin{proposition}\label{prop3_7}
{\bf [Monotony and injectivity]}
Let $V,W$ be SI subspaces.
\par
(i) $V\subset W$ iff $\tau_{V,T}\leq\tau_{W,T}$ a.e. for all positive operators $T$ iff
$\tau_{V,f}\leq\tau_{W,f}$ a.e.
for all $f\in\ltwozn$.
\par
(ii) $V=W$ iff $\tau_{V,T}=\tau_{W,T}$ a.e. for all positive operators $T$ iff $\tau_{V,f}=\tau_{W,f}$ a.e for
all $f\in\ltwozn$.
\end{proposition}

\begin{proof}
(i) It is clear that the first statement implies the second; just use the additivity property for the SI spaces
$V$ and $W\ominus V$. The third statement is just a particular case of the second one. So the only interesting implication is
from the third statement to the first one.
\par
Let $J_V$ and $J_W$ be the corresponding periodic range functions for $V$ and $W$. The hypothesis implies, according to
lemma \ref{lem21_9_1}, that
$$J_V(\xi)\subset J_W(\xi),\quad(\mbox{ for a.e. }\xi\in\mathbb{R}^n),$$
and this implies in turn that $V\subset W$ (just look at theorem \ref{th2_1}).
\par
(ii) is a consequence of (i) by a double inclusion argument.
\end{proof}

\par
For $a\in\mathbb{R}^n$ we define the modulation of $f\in\ltworn$ by
$$M_a(f)(x)=e^{i\ip{a}{x}}f(x),\quad(x\in\mathbb{R}^n).$$
The local trace function behaves nicely under modulation. This is expressed in the next proposition.

\begin{proposition}\label{prop3_8}
{\bf [Modulation]}
Let $V$ be a SI subspace and $a\in\mathbb{R}^n$. Then $M_aV$ is a SI subspace and for
all positive operators $T$ on $\ltwozn$ and
all vectors $f\in\ltwozn$:
\begin{equation}\label{eq3_8_1}
\tau_{M_aV,T}(\xi)=\tau_{V,T}(\xi-a),\quad\mbox{ for a.e. }\xi\in\mathbb{R}^n;
\end{equation}
\begin{equation}\label{eq3_8_2}
\tau_{M_aV,f}(\xi)=\tau_{V,f}(\xi-a),\quad\mbox{ for a.e. }\xi\in\mathbb{R}^n.
\end{equation}
\end{proposition}

\begin{proof}
The modulation and the translations satisfy a commutation relation: for $k\in\mathbb{Z}^n$ and $\varphi\in\ltworn$,
$$T_kM_a\varphi(x)=e^{i\ip{a}{x-k}}\varphi(x-k)=e^{-i\ip{a}{k}}M_aT_k\varphi(x),\quad(x\in\mathbb{R}^n).$$
This relation shows that $M_aV$ is shift invariant.
\par
Now take a NTF generator $\phi$ for $V$ (it exists by theorem \ref{th1_3}). Then, as $M_a$ is unitary,
$$\{M_aT_k\varphi\,|\,\varphi\in\phi\}$$
is a NTF for $M_aV$. Using the commutation relation and the fact that $e^{-i\ip{a}{k}}$ are just constants of
modulus 1, we
see that
$$\{T_kM_a\varphi\,|\,\varphi\in\phi\}$$
is a NTF for $M_aV$. Therefore we can safely use theorem \ref{th3_3} and compute:
$$\tau_{M_aV,T}(\xi)=\sum_{\varphi\in\phi}\ip{T\Tper(M_a\varphi)(\xi)}{\Tper(M_a\varphi)(\xi)}.$$
But
$$\Tper
M_a\varphi(\xi)=\left(\widehat{M_a\varphi}(\xi+2k\pi)\right)_{k\in\mathbb{Z}^n}=\left(\widehat{\varphi}(\xi-a+2k\pi)\right)_{k\in\mathbb{Z}^n}$$
$$=\Tper\varphi(\xi-a),$$
and (\ref{eq3_8_1}) follows with (\ref{eq3_8_2}) as its consequence.
\end{proof}

\par
The dilation by an $n\times n$ non-singular matrix $A$ is the unitary operator on $\ltwor$ defined by
$$D_Af(x)=|\operatorname*{det}A|^{\frac{1}{2}}f(Ax),\quad(x\in\mathbb{R}^n,f\in\ltworn).$$
We will consider only matrices $A$ which preserve the lattice $\mathbb{Z}^n$, because in this case $D_AV$ is
shift invariant whenever $V$ is.

\begin{proposition}\label{prop3_9}
{\bf [Dilation]}
Let $V$ be a SI subspace and $A$ an $n\times n$ integer matrix with $\operatorname*{det}A\neq 0$. Then $D_AV$ is
shift invariant and,
for every positive operator $T$ on $\ltwozn$ and every vector $f\in\ltwozn$:
\begin{equation}\label{eq3_9_1}
\tau_{D_AV,T}(\xi)=\sum_{d\in\mathcal{D}}\tau_{V,D_d^*TD_d}\left(\left(A^*\right)^{-1}(\xi+2d\pi)\right),\quad\mbox{ for a.e.
}\xi\in\mathbb{R}^n,
\end{equation}
\begin{equation}\label{eq3_9_2}
\tau_{D_AV,f}(\xi)=\sum_{d\in\mathcal{D}}\tau_{V,D_d^*f}\left(\left(A^*\right)^{-1}(\xi+2d\pi)\right),\quad\mbox{
for
a.e. }\xi\in\mathbb{R}^n,
\end{equation}
where $\mathcal{D}$ is a complete set of representatives of the cosets $\mathbb{Z}^n/A^*\mathbb{Z}^n$
and
$D_d$ is the linear operator on $\ltwozn$ defined by
$$(D_d\alpha)(k)=\left\{\begin{array}{ccc}
\alpha(l),&\mbox{if}&k=d+A^*l\\
0,& & otherwise
\end{array}
\right.,\quad(k\in\mathbb{Z}^n,\alpha\in\ltwozn).$$
\end{proposition}

\begin{proof}
The dilation and the translation satisfy the following commutation relation which can be easily verified:
$$T_kD_A=D_AT_{Ak},\quad(k\in\mathbb{Z}^n).$$
This shows that $D_AV$ is shift invariant.
\par
We can decompose $V$ as the orthogonal sum $V=\oplus_{i\in I}S(\varphi_i)$, where $\varphi_i$ is a quasi-orthogonal generator
of $S(\varphi_i)$ (see theorem \ref{th1_3}). Since $D_A$ is unitary, $D_AV=\oplus_{i\in I}D_AS(\varphi_i)$ and,
using the
additivity property of the local trace function (proposition \ref{prop3_6}), it suffices to prove the formula (\ref{eq3_9_1})
just for the case when $V=S(\varphi)$ with $\varphi$ quasi-orthogonal generator for $S(\varphi)$. We will
assume this is the case.
\par
From the commutation relation we see that, if $\mathcal{L}$ is a
complete set of $|\operatorname*{det}A|$ representatives of the
cosets $\mathbb{Z}^n/A\mathbb{Z}^n$, then
$$\{D_AT_l\varphi\,|\, l\in\mathcal{L}\}$$
will span $D_AV$ by translations.
\par
For $l\in\mathcal{L}$, consider
\begin{align*}
\phi_l(\xi)&=\Tper(D_AT_l\varphi)(\xi)\\
&=\left(|\operatorname*{det}
A|^{-\frac{1}{2}}\widehat{\varphi}\left(\left(A^*\right)^{-1}(\xi+2k\pi)\right)
e^{-i\ip{(A^*)^{-1}(\xi+2k\pi)\,}{\,l}}\right)_{k\in\mathbb{Z}^n}.
\end{align*}
For $d\in\mathcal{D}$ define $\psi_d\in L^2_{per}(\Rn,\ltwozn)$ by
\begin{align*}
\psi_d(\xi)(k)&=\left\{\begin{array}{ccc}
\widehat{\varphi}((A^*)^{-1}(\xi+2k\pi)),&\mbox{if}&k\in d+A^*\Zn\\
0,& &\mbox{otherwise}
\end{array}
\right.\\
&=\left\{\begin{array}{ccc}
\widehat{\varphi}((A^*)^{-1}(\xi+2d\pi)+2\pi l),&\mbox{if}&k= d+A^*l,\mbox{ with }l\in\Zn\\
0,& &\mbox{otherwise}
\end{array}
\right.\\
&=D_d\left(\Tper\varphi((A^*)^{-1}(\xi+2\pi d))\right).
\end{align*}
Then $\phi_l(\xi)$ and $\psi_d(\xi)$ are related by the following linear equations:
$$\phi_l(\xi)=e^{-i\ip{(A^*)^{-1}\xi\,}{\,l}}|\operatorname*{det}A|^{-1/2}
\sum_{d\in\mathcal{D}}e^{-i\ip{(A^*)^{-1}2\pi d\,}{\,l}}\psi_d(\xi).$$
Since the $|\operatorname*{det}A|\times |\operatorname*{det}A|$ matrix
$$\left( |\operatorname{det}A|^{-1/2}e^{-i\ip{(A^*)^{-1}2\pi
d\,}{\,l}}\right)_{d\in\mathcal{D},l\in\mathcal{L}}$$
is unitary, it follows that $\{\psi_d(\xi)\,|\, d\in\mathcal{D}\}$ and $\{\phi_l(\xi)\,|\, l\in\mathcal{L}\}$ span the same
subspace of $\ltwozn$, namely $J_{D_AV}(\xi)$ (use theorem \ref{th2_1}).
\par
$\varphi$ is a quasi-orthogonal generator so $\operatorname*{Per}|\widehat{\varphi}|^2$ is a characteristic
function that is
$\|\Tper\varphi(\xi)\|_{\ltwozn}$ is either 0 or 1, and as $D_d$ is an isometry, $\|\psi_d(\xi)\|_{\ltwozn}\in\{0,1\}$.
Also $\psi_d(\xi)$ and $\psi_{d'}(\xi)$ are perpendicular when $d\neq d'$ and, in conclusion $\{\psi_d(\xi)\,|\, d\in\mathcal{D}\}$ is
a NTF for $J_{D_AV}(\xi)$. Therefore we can use these vectors to compute the local trace function (proposition
\ref{prop21_4}):
\begin{align*}
\tau_{D_AV,T}(\xi)&=\sum_{d\in\mathcal{D}}\ip{T\psi_d(\xi)}{\psi_d(\xi)}\\
&=\sum_{d\in\mathcal{D}}\ip{TD_d\left(\Tper\varphi((A^*)^{-1}(\xi+2\pi
d))\right)}{D_d\left(\Tper\varphi((A^*)^{-1}(\xi+2\pi d))\right)}\\
&=\sum_{d\in\mathcal{D}}\ip{D_d^*TD_d\left(\Tper\varphi((A^*)^{-1}(\xi+2\pi
d))\right)}{\Tper \varphi((A^*)^{-1}(\xi+2\pi d))}\\
&=\sum_{d\in\mathcal{D}}\tau_{V,D_d^*TD_d}((A^*)^{-1}(\xi+2\pi d))
\end{align*}
This proves (\ref{eq3_9_1}). (\ref{eq3_9_2}) follows from (\ref{eq3_9_1}) because $D_d^*P_fD_d=P_{D_d^*f}$.
\end{proof}
\par
As we promised, the local trace function incorporates the dimension function and the spectral function defined by
M.Bownik and Z.Rzeszotnik in \cite{BoRz}.
\par
Recall that for a shift invariant subspace $V$ of $\ltworn$, its dimension function is defined as
$${\operatorname*{dim}}_V(\xi)=\operatorname*{dim}J_{per}(\xi),\quad(\xi\in\Rn),$$
where $J_{per}$ is the periodic function associated to $V$.
\par
The spectral function introduced in \cite{BoRz} is defined by
$$\sigma_V(\xi+2k\pi)=\|J_{per}(\xi)\delta_k\|^2,\quad(\xi\in[-\pi,\pi]^n,k\in\Zn),$$
where $\delta_k\in\ltwozn$, $\delta_k(l)=\left\{\begin{array}{ccc}
1,&\mbox{if}&k=l\\
0,& &\mbox{otherwise}
\end{array}
\right.$
$(k\in\mathbb{Z}^n)$.
\par
For a similar treatment of the dimension function, using the
Gramian the reader can also consult \cite{RS4}.

\begin{proposition}\label{prop3_10}
Let $V$ be a SI space. \\
(i)
\begin{equation}\label{eq3_10_1}
\tau_{V,I}={\operatorname*{dim}}_V.
\end{equation}
(ii)
\begin{equation}\label{eq3_10_2}
\tau_{V,\delta_0}=\sigma_V.
\end{equation}
\end{proposition}

\begin{proof}
Let $J_{per}$ be the periodic range function of $V$.
$$\tau_{V,I}(\xi)=\Trace(IJ_{per}(\xi))=\operatorname*{dim}J_{per}(\xi)=\dim_V(\xi).$$
\par
For $\xi\in [-\pi,\pi)^n$, by proposition \ref{prop3_1},
$$\tau_{V,\delta_0}(\xi)=\|J_{per}(\xi)\delta_0\|^2=\sigma_V(\xi).$$
If $k\in\mathbb{Z}^n$ then, using the periodicity of the local trace function stated in
proposition \ref{prop3_4},
$$\tau_{V,\delta_0}(\xi+2k\pi)=\tau_{V,\lambda(k)\delta_0}(\xi)=\tau_{V,\delta_k}(\xi)=
\|J_{per}(\xi)\delta_k\|^2=\sigma_V(\xi).$$
\end{proof}

\section{\label{convergence}Convergence theorems}
In this section we study the behavior of the local trace function with respect to limits.
More precisely, we consider the following question: if $(V_j)_{j\in\mathbb{N}}$ is a sequence
of SI subspaces of $\ltworn$ such that $P_{V_j}$ converges in the strong operator topology to
$P_V$ for some SI subspace $V$, what can be said about the convergence of the local trace
functions $\tau_{V_i,T}$?
\par
We begin with a lemma and then give some useful partial answers to this question.
\begin{lemma}\label{le4_1}
Suppose $(V_j)_{j\in\mathbb{N}}$ and $V$ are SI subspaces of $\ltworn$ such that
$P_{V_j}$ converges in the strong operator topology to $P_V$. Then for any $f\in\ltworn$ and
any $A>0$
$$\int_{[-A,A]^n}\|J_{V_j}(\xi)(\Tper f(\xi))-J_V(\xi)(\Tper f(\xi))\|^2\,d\xi\mbox{ converges to
}0\mbox{ as }j\rightarrow\infty,$$
where $J_{V_j}$ and $J_{V}$ are the corresponding periodic range functions.
\end{lemma}
\begin{proof}
We can consider $A=\pi$ because then the result is obtained using the periodicity of the range
function. With proposition \ref{prop2_4},
$$J_{V_j}(\xi)(\Tper f(\xi))=\Tper(P_{V_j}f)(\xi),\quad (\xi\in\mathbb{R}^n).$$
And then
$$\int_{[-\pi,\pi]^n}\|J_{V_j}(\xi)(\Tper f(\xi))-J_V(\xi)(\Tper f(\xi))\|^2\,d\xi$$
$$=\int_{[-\pi,\pi]^n}\|\Tper(P_{V_j}f)(\xi)-\Tper(P_Vf)(\xi)\|^2\,d\xi$$
$$=(2\pi)^n\|P_{V_j}f-P_Vf\|^2_{\ltworn}\rightarrow 0,\mbox{ as }j\rightarrow\infty.$$
\end{proof}

\begin{theorem}\label{th4_2}
Let $(V_j)_{j\in\mathbb{N}}$ and $V$ be SI subspaces such that $P_{V_j}$ converges to $P_V$ in
the strong operator topology. Then, for any bounded measurable set $E$ and every $f\in\ltwozn$
$$\int_E|\tau_{V_j,f}(\xi)-\tau_{V,f}(\xi)|\,d\xi\rightarrow 0,\mbox{ as
}j\rightarrow\infty.$$
\end{theorem}

\begin{proof}
Let $J_{V_j}$ and $J_V$ be the corresponding periodic range function. We can assume
$E=[-\pi,\pi]^n$ because then the result is a consequence of the periodicity of the local
trace function (proposition \ref{prop3_4}) and we can assume also $\|f\|=1$. Using lemma \ref{le4_1} for a
$g\in\ltworn$
with
$\mathcal{T}g(\xi)=f$ for all $\xi\in[-\pi,\pi]^n$, we have
\begin{align*}
0&=\lim_{j\rightarrow\infty}\int_{[-\pi,\pi]^n}\|J_{V_j}(\xi)f-J_V(\xi)f\|^2\,d\xi\\
&\geq\limsup_{j\rightarrow\infty}\int_{[-\pi,\pi]^n}\left|\|J_{V_j}(\xi)f\|-\|J_V(\xi)f\|\right|^2\,d\xi\\
&=\limsup_{j\rightarrow\infty}\int_{[-\pi,\pi]^n}|\tau_{V_j,f}^{1/2}(\xi)-\tau_{V,f}^{1/2}(\xi)|^2\,d\xi\\
&\geq\frac{1}{4}\limsup_{j\rightarrow\infty}\int_{[-\pi,\pi]^n}|\tau_{V_j,f}(\xi)-\tau_{V,f}(\xi)|^2\,d\xi\\
&\geq\frac{1}{4}(2\pi)^n\limsup_{j\rightarrow\infty}\int_{[-\pi,\pi]^n}|\tau_{V_j,f}(\xi)-\tau_{V,f}(\xi)|\,d\xi,
\end{align*}
where for the third line we used proposition \ref{prop3_1}, for the fourth we used the fact
that $\tau_{V,f},\tau_{V_j,f}$ are less then $\|f\|^2=1$ (proposition \ref{prop3_5} (i)) and for the last one we
used Holder's inequality.
\end{proof}

\begin{theorem}\label{th4_3}
Let $(V_j)_{j\in\mathbb{N}}$ be an increasing sequence of SI subspaces of $\ltworn$,
$$V:=\overline{\bigcup_{j\in\mathbb{N}}V_j},$$
$T$ a positive operator on $\ltwozn$ and $f\in\ltwozn$. Then, for a.e. $\xi\in\mathbb{R}^n$,
$\tau_{V_j,T}(\xi)$ increases to $\tau_{V,T}(\xi)$ and $\tau_{V_j,f}(\xi)$ increases to
$\tau_{V,f}(\xi)$.
\end{theorem}

\begin{proof} The monotony is a consequence of proposition \ref{prop3_7}.
\par
Again, we denote by $J_{V_j}$ and $J_V$ the corresponding periodic range function. By
proposition \ref{prop2_3}, $(J_{V_j}(\xi))_{j\in\mathbb{N}}$ is increasing and
$$J_V(\xi)=\overline{\bigcup_{j\in\mathbb{N}}J_{V_j}(\xi)},$$
for almost every $\xi\in\mathbb{R}^n$.
\par
Apply now proposition \ref{prop21_6} to our situation to conclude that
$$\tau_{V,T}(\xi)=\Trace(TJ_V(\xi))=\lim_{j\rightarrow\infty}\Trace(TJ_{V_j}(\xi))=\lim_{j\rightarrow\infty}\tau_{V_j,T}(\xi).$$
The second statement is a particular case of the first.
\end{proof}

\section{\label{wavelets}Wavelets and the local trace function}
We reserved this section for applications of the local trace function to wavelets.
Throughout this section $A$ will be an $n\times n$ dilation matrix (i.e all eigenvalues $\lambda$ have
$|\lambda|>1$), and $A$ preserves the lattice $\Zn$ (i.e. $A\Zn\subset \Zn$).
\begin{remark}\label{rem5_2}
The local trace function can be used to obtain the characterization of wavelets (see \cite{Bo2} or \cite{Ca}), namely:
\par
A set $\Psi=\{\psi^1,...,\psi^L\}$ is a NTF wavelet (i.e. the set
$$X(\Psi):=\{D_A^jT_k\psi^l\,|\,j\in\mathbb{Z},k\in\Zn,l\in\{1,...,L\}\}$$
is a NTF for $\ltworn$) if and only if the following equations hold for a.e. $\xi\in\Rn$:
\begin{equation}\label{eq5_2_1}
\sum_{\psi\in\Psi}\sum_{j=-\infty}^\infty|\widehat{\psi}|^2((A^*)^j\xi)=1,
\end{equation}
\begin{equation}\label{eq5_2_2}
\sum_{\psi\in\Psi}\sum_{j\geq0}\widehat{\psi}((A^*)^j\xi)\overline{\widehat{\psi}}((A^*)^j(\xi+2s\pi))=0,\quad(s\in\Zn\setminus A^*\Zn).
\end{equation}

 The argument used here will be the same as the one used in \cite{Bo2}, the only difference is that, instead of the
Gramian we employ the local trace function. Here is a sketch of
the proof: It is known (see \cite{CSS}) that $X(\Psi)$ is an
affine NTF for $\ltworn$ if and only if the quasi-affine system
$X^q(\Psi)$ is a NTF for $\ltworn$. Recall that
$$X^q(\Psi):=\{\tilde{\psi}^l_{j,k}\,|j\in\mathbb{Z},k\in\Zn,l\in\{1,...,L\}\},$$
with the convention
$$\tilde{\psi}_{j,k}(x)=\left\{\begin{array}{ccc}
D_A^jT_k\psi(x),&\mbox{if}&j\geq0,k\in\Zn\\
|\operatorname*{det} A|^{j/2}T_kD_A^j\psi(x)&\mbox{if}&j<0,k\in\Zn.
\end{array}\right.$$
This can be reformulated as
$$\{\tilde{\psi}^l_{j,0}\,|\, j<0\}\cup\{\tilde{\psi}_{j,r}^l\,|\,j\geq0,r\in\mathcal{L}_j\}$$
is a NTF generator for $\ltworn$ (here $\mathcal{L}_j$ is a complete set of representatives of $\Zn/A^j\Zn$.
\par
Now use corollary \ref{cor3_3_1}, and all one has to do is a computation. The computation is done in lemma 2.3 in
\cite{Bo2} and one obtains, for a.e. $\xi\in\Rn$,
$$\sum_{\psi\in\Psi}\sum_{j=-\infty}^\infty|\widehat{\psi}|^2((A^*)^j(\xi)=1,$$
$$t_{(A^*)^{-m}p}((A^*)^{-m}(\xi))=0,\quad(p\in\Zn,p\neq0),$$
where
$$t_s(\xi):=\sum_{\psi\in\Psi}\sum_{j\geq0}\widehat{\psi}((A^*)^j\xi)\overline{\widehat{\psi}}((A^*)^j(\xi+2s\pi)),\quad(\xi\in\Rn,s\in\Zn\setminus
A^*\Zn)$$
and $m=\max\{j\in\mathbb{Z}\,|\,(A^*)^{-j}p\in\Zn\}$. This leads immediately to the equivalence to (\ref{eq5_2_1}) and
(\ref{eq5_2_2}).
\par
Another observation is that, if we know that the equivalence "
$X(\psi)$ is a NTF for $\ltworn$ iff (\ref{eq5_2_1}) and
(\ref{eq5_2_2}) hold", then the argument above also shows that the
affine system $X(\psi)$ is a NTF for $\ltworn$ iff $X^q(\psi)$ is
a NTF for $\ltworn$.
\end{remark}

\begin{theorem}\label{th5_1}
Let $\Psi=\{\psi^1,...,\psi^L\}$ be a semiorthogonal wavelet, i.e. the affine system
$$\{D_{A}^jT_k\psi\,|\, j\in\mathbb{Z}, k\in\mathbb{Z}^n,\psi\in\Psi\}$$
is a NTF for $\ltworn$ and $W_i\perp W_j$ for $i\neq j$, where
$$W_j=\overline{\operatorname*{span}}\{D_{A}^jT_k\psi\,|\,k\in\mathbb{Z}^n,\psi\in\Psi\}=D_{A}^j(S(\Psi)),\quad(j\in\mathbb{Z}),$$
$$V_j=\bigoplus_{i<j}W_i.$$
Then \\
(i) The set
$$\{|\operatorname*{det}A|^{j/2}T_kD_{A}^j\psi\,|\,j<0,k\in\mathbb{Z}^n\}$$
is a NTF for $V_0$.\\
(ii) For $\psi\in\Psi$, $j\geq 1$ and $\xi\in\mathbb{R}^n$, denote by
$f_{\psi}^j(\xi)$ the vector in $\ltwozn$ defined by
$$f_{\psi}^j(\xi)(k)=\widehat{\psi}((A^*)^j(\xi+2k\pi)),\quad(k\in\Zn).$$
Then, for almost every $\xi\in\mathbb{R}^n$, the set
$$\{f_{\psi}^j(\xi)\,|\, j\geq 1,\psi\in\Psi\}$$
is a NTF for $J_{V_0}(\xi)$, where $J_{V_0}$ is the periodic range function of $V_0$.\\
(iii) For every positive operator $T$ on $\ltwozn$, every $f\in\ltwozn$ and almost every
$\xi\in\mathbb{R}^n$,
\begin{equation}\label{eq5_1_1}
\tau_{V_0,f}(\xi)=\sum_{\psi\in\Psi}\sum_{j=1}^\infty\ip{Tf_{\psi}^j(\xi)}{f_{\psi}^j(\xi)},
\end{equation}
\begin{equation}\label{eq5_1_2}
\tau_{V_0,f}(\xi)=\sum_{\psi\in\Psi}\sum_{j=1}^\infty|\sum_{k\in\mathbb{Z}^n}\widehat{\psi}((A^*)^j(\xi+2k\pi))\overline{f}_k|^2.
\end{equation}
\end{theorem}

\begin{proof}
(i) follows, as we mentioned before, from the equivalence between
affine frames and quasi-affine frames (see \cite{CSS} or theorem 1.4 in \cite{Bo2}) and the
orthogonality relations given in the hypothesis.
\par
(ii) If we compute the Fourier transform of $|\operatorname*{det}A|^{-j/2}D_{A}^{-j}\psi$ for $j>0$ we get
$\widehat{\psi}((A^*)^j\xi)$ and (ii) is obtained from (i) and theorem \ref{th2_5}.
\par
(iii) The fact asserted in (ii) shows that we can use the vectors $f_{\psi}^j$ to compute the
trace and the resulting formulas are exactly (\ref{eq5_1_1}) and (\ref{eq5_1_2}).
\end{proof}

\begin{remark}\label{rem5_3}
Another fundamental fact from the theory of wavelets can be obtained with the aid of the local trace function:
the dimension function equals the trace function (see \cite{Web}). Let's recall some notions. Consider
$\Psi:=\{\psi^1,...,\psi^L\}$ a semi-orthogonal wavelet (as in theorem \ref{th5_1}). We keep the notations
for $V_j$ and $W_j$. The dimension function associated to $\Psi$ is
$$D_{\Psi}(\xi)=\sum_{k\in\Zn}\sum_{\psi\in\Psi}\sum_{j\geq
1}|\widehat{\psi}|^2\left((A^*)^j(\xi+2k\pi)\right).$$
\par
The definition of the multiplicity function requires a little bit of harmonic analysis (for details look in
\cite{BMM} or \cite{B}). The subspace $V_0$ in
invariant under translations by integers (shift invariant). So one has a unitary representation of the
locally compact abelian group $\Zn$ on $V_0$ by translations. The Stone-Mackey theory shows that this
representation is determined by a projection-valued measure which in turn is determined by a positive measure
class on the dual group $\widehat{\Zn}$ and a measurable multiplicity function
$m_{V_0}:\widehat{\Zn}\rightarrow\{0,1...,\infty\}$.
$\widehat{\Zn}$ can be identified with $[-\pi,\pi]^n$. In our case, the measure is the Lebesgue measure, so
the determinant is the multiplicity function $m_{V_0}$. The beautiful result is
\begin{equation}\label{eq5_3_1}
D_{\Psi}(\xi)=m_{V_0}(\xi),\quad(\xi\in[-\pi,\pi]^n).
\end{equation}
We use the local trace function to prove it. First, if we use the equation (\ref{eq5_1_1}) with $T=I$, the
identity on $\ltwozn$, we obtain
\begin{equation}\label{eq5_3_2}
\dim_{V_0}(\xi)=\tau_{V_0,I}(\xi)=D_{\Psi}(\xi),\quad(\xi\in[-\pi,\pi]^n).
\end{equation}
To equate the local trace function with the multiplicity function we use a result from \cite{BM}: if
$$S_j:=\{\xi\in[-\pi,\pi]^n\,|\,m_{V_0}(\xi)\geq j\},\quad(j\in\mathbb{N},j\geq 1)$$
then there exist a NTF generator $\{\varphi_j\,|\,j\geq 1\}$ for $V_0$ such that
\begin{equation}\label{eq5_3_2_1}
\sum_{k\in\Zn}\widehat{\varphi}_i\overline{\widehat{\varphi}}_j(\xi+2k\pi)=
\left\{\begin{array}{ccc}
\chi_{S_i}(\xi),&\mbox{if}&i=j\\
0,&\mbox{if}&i\neq j.
\end{array}
\right.
\end{equation}
We can use theorem \ref{th3_3} for $T=I$ and we have for a.e. $\xi$:
$$\tau_{V_0,I}(\xi)=\sum_{j\geq 1}\ip{\Tper\varphi_j(\xi)}{\Tper\varphi_j(\xi)
}=\sum_{j\geq
1}\sum_{k\in\Zn}|\widehat{\varphi}_j|^2(\xi)=\sum_{j\geq 1}\chi_{S_j}(\xi).$$
Therefore
\begin{equation}\label{eq5_3_3}
\dim_{V_0}(\xi)=\tau_{V_0,I}(\xi)=m_{V_0}(\xi),\quad(\xi\in[-\pi,\pi]^n).
\end{equation}
Consequently the multiplicity function and the dimension function are just two disguises of the local trace
function at $T=I$.
\end{remark}
The next theorem gives the equations that relates scaling
functions (i.e NTF generators for $V_0$) wavelets.
\begin{theorem}\label{th5_4}
Let $\Psi=\{\psi^1,...,\psi^L\}$ be a semi-orthogonal wavelet as in theorem \ref{th5_1}. Let $\Phi$ be a countable subset
of $\ltworn$. The following affirmations are equivalent:
\begin{enumerate}
\item
$\Phi$ is contained in $V_0$ and is a NTF generator for $V_0$;
\item
The following equations hold: for every $s\in\Zn$,
\begin{equation}\label{eq5_4_2}
\sum_{\psi\in\Psi}\sum_{j\geq 1}\widehat{\psi}((A^*)^j\xi)\overline{\widehat{\psi}}((A^*)^j(\xi+2s\pi))=
\sum_{\varphi\in\Phi}\widehat{\varphi}(\xi)\overline{\widehat{\varphi}}(\xi+2s\pi).
\end{equation}
for a.e. $\xi\in\Rn$.
\end{enumerate}
\end{theorem}

\begin{proof}
Use theorem \ref{th3_3_2} and theorem \ref{th5_1}.
\end{proof}

\begin{corollary}\label{cor5_5}
Let $\Psi=\{\psi^1,...,\psi^L\}$ be a semi-orthogonal wavelet as in theorem \ref{th5_1} and
$\Phi$ a NTF generator for $V_0$. Then
\begin{equation}\label{eq5_5_1}
\sum_{\psi\in\Psi}|\widehat{\psi}|^2(\xi)=\sum_{\varphi\in\Phi}|\widehat{\varphi}|^2((A^*)^{-1}\xi)-
\sum_{\varphi\in\Phi}|\widehat{\varphi}|^2(\xi),\quad\mbox{a.e. on }\Rn;
\end{equation}
For all $s\in\Zn\setminus A^*\Zn$,
\begin{equation}\label{eq5_5_2}
\sum_{\psi\in\Psi}\widehat{\psi}(\xi)\overline{\widehat{\psi}}(\xi+2s\pi)=
-\sum_{\varphi\in\Phi}\widehat{\varphi}(\xi)\overline{\widehat{\varphi}}(\xi+2s\pi),\quad\mbox{a.e. on }\Rn.
\end{equation}
\end{corollary}
\begin{proof}
For (\ref{eq5_5_1}), take $s=0$, write the equation (\ref{eq5_4_2}) for $\xi$ and $(A^*)^{-1}\xi$, and substract the first from
the second. For (\ref{eq5_5_2}), substract equation (\ref{eq5_4_2}) from equation (\ref{eq5_2_2}).
\end{proof}

\end{document}